\documentclass[12pt]{article}
\newcommand{\copyleft}{
GNU FDL\thanks{
Copyright (C) 1994, 1998, 2009 Peter G. Doyle.
Permission is granted to copy, distribute and/or modify this document
under the terms of the GNU Free Documentation License, 
as published by the Free Software Foundation;
with no Invariant Sections, no Front-Cover Texts, and no Back-Cover Texts.
}}
\title{Frustration solitaire}
\author{Peter G. Doyle \and Charles M. Grinstead \and J. Laurie Snell}
\date{Version dated 2 April 2009
\\ \copyleft
}
\begin{document}
\maketitle

\begin{abstract}
In this expository article,
we discuss the 
rank-derangement problem, which asks for the number of permutations
of a deck of cards such that each card is replaced by a card of a
different rank.
This combinatorial problem arises in computing the probability
of winning the game of `frustration solitaire'.
The solution is a prime example of the method of inclusion and exclusion.
We also discuss and announce the solution to 
Montmort's `Probleme du Treize', a related problem dating back to circa 1708.

This revised version incorporates corrections requested by
Steven Langfelder to the historical remarks.
\end{abstract}

\section{Frustration solitaire and rank-derangements}

In the game of frustration solitaire,
you shuffle a deck of cards thoroughly--at least 13 times--so that by the
time you're done all orderings of the deck are equally likely.
Now you run through the deck, turning over the cards
one at a time as you call out, `Ace, two, three, four, five, six, seven, eight,
nine, ten, jack, queen, king, ace, two, three, \ldots ,'
and so on, so that you end up calling out the thirteen ranks four times each.
If the card that comes up ever matches the rank you call out as you turn it
over, then you lose.

The
{\em rank-derangement problem} asks for the
number of permutations of a deck of cards such that after the
permutation, every card has been
replaced by a card of a different rank.
We denote this number by $R_{13}$, where the subscript $13$ is there
to remind us that there are 13 different ranks.
More generally, the rank-derangement problem asks for the $R_n$,
the number of rank-derangements of a $4n$-card deck
containing cards of $n$ different ranks in each of four suits.
Of course we could generalize the problem still
further by varying the number of suits,
but we won't get into that here.

To see the connection between the rank-derangement problem  and
frustration solitaire, imagine that the deck
starts out in the order: ace through king of spades, ace through king of hearts,
ace through king of diamonds, ace through king of clubs.
Starting with this ordering, you win if the permutation you do when you
shuffle is a rank-derangement,
which happens with probability
$R_{13}/52!$.

Now of course you probably didn't start with this particular ordering,
unless you just won at one of the more conventional kinds of solitaire.
But because we're assuming that after the shuffle all orderings are
equally likely, your probability of winning doesn't depend on
what order the deck started in, so
your probability of winning is still going to be
$R_{13}/52!$.

By the way, the special ordering we've just been discussing is not at
all the ordering of a brand new deck of Bees or Bicycles,
or practically any other brand of top-quality playing cards.
New decks come in quite a different ordering, and the particular ordering
that they come in is
very important if you are going to be playing certain special games
like {\em bore} or {\em new age solitaire}, and are only planning
on shuffling the deck seven or eight times.
We're not going to say anything more about that here.
We just didn't want to give you the wrong impression about the order of the
cards in a new deck.

\section{Historical background}
The roots of the rank-derangement go back to a gambling game closely
resembling frustration solitaire that was studied long ago by the
Chevalier de Montmort.

In 1708 Pierre de Montmort published 
the first edition of his book {\em Essay d'Analyse sur les Jeus de
Hazard}
\cite{montmort:hazard}
(Analytical Essay on Games of Chance).  He was inspired by
recent work of James Bernoulli in probability.  Montmort
hoped that, by applying techniques developed by Bernoulli to analyze 
common card and dice games, he could show people
that certain of their methods of play based on superstitions should 
replaced by more rational behavior. One of the games that he analyzed 
was {\em Jeu de Treize} (Thirteen).  This game  was played as follows:  

One person is chosen as  dealer and the others are players.  Each
player puts up a stake. The dealer shuffles the cards and turns them
up one at a time calling out, `Ace, two, three,\ldots, king', just as in 
frustration solitaire. If the dealer goes through the 13 cards 
without a match he (or she?)
pays the players an amount equal to their stake, and the deal passes
to someone else.  If there is a match the dealer collects the
players' stakes; the players put up new stakes, and the dealer continues
through the deck, calling out, `Ace, two, three, \ldots'.
If the dealer runs out of cards he reshuffles and continues the count where he
left off. He continues until there is a run of 13 without a match and then a 
new dealer is chosen.

Montmort's `Probl\`eme du Treize' was to find the expected value of this
game to the dealer.
The answer, found by Peter Doyle in 1994,
\begin{verbatim}
26516072156010218582227607912734182784642120482136\
09144671537196208993152311343541724554334912870541\
44029923925160769411350008077591781851201382176876\
65356317385287455585936725463200947740372739557280\
74593843427478766496507606399053826118938814351354\
73663160170049455072017642788283066011710795363314\
27343824779227098352817532990359885814136883676558\
33113244761533107206274741697193018066491526987040\
84383914217907906954976036285282115901403162021206\
01549126920880824913325553882692055427830810368578\
18861208758248800680978640438118582834877542560955\
55066287892712304826997601700116233592793308297533\
64219350507454026892568319388782130144270519791882   /
33036929133582592220117220713156071114975101149831\
06336407213896987800799647204708825303387525892236\
58132301562800562114342729062565897443397165719454\
12290800708628984130608756130281899116735786362375\
60671849864913535355362219744889022326710115880101\
62859313519792943872232770333969677979706993347580\
24236769498736616051840314775615603933802570709707\
11959696412682424550133198797470546935178093837505\
93488858698672364846950539888686285826099055862710\
01318150621134407056983214740221851567706672080945\
86589378459432799868706334161812988630496327287254\
81845887935302449800322425586446741048147720934108\
06135061350385697304897121306393704051559533731591
\end{verbatim}
which as you can see is something like .8. 
The solution, which we will not discuss here, involves
having the computer solve something like 4 million
variations on the frustration solitaire problem.

Montmort first considered  the problem of finding the
probability that there would be a match before getting
through the first 13 cards. He started by assuming
that the deck had only 13 cards all of the same suit and
showed
that the probability of getting no match is
very close to $1/e= .3678\ldots$,
so the probability of getting a match is very close to .6321\ldots.
He used the method of recursion to derive this result, thus giving
the first solution to the problem we now call the `derangement problem',
or `hat-check problem',
which we will discuss below.

In later correspondence with Nicholas Bernoulli, Montmort found that,
for a normal 52 card deck, the probability of getting through 13 cards without
a match is .356\ldots,
so the dealer wins the first round of Treize with
probability .643\ldots. 
(Note that the dealer does a little better on the first round with a
52 card deck than with a 13 card deck.)
Thus Montmort showed that the dealer has a
significant advantage even without considering the additional winnings
from further rounds before giving up the deal.

The game of Treize is still mentioned in some books on card games,
but in forms not so advantageous to the dealer.
In one  version, the dealer deals  until there is
a match or until 13 cards have come up without a match.  If there is a match 
on card $n$, he or she wins \$$n$ from each player.

The current interest in the problem came from a column of Marilyn vos
Savant
\cite{savant:frustration}.
Charles Price wrote to ask about his experience
playing 
frustration solitaire.   He found that he
he  rarely won and wondered how often he
should win. Marilyn answered by  remarking that the expected number of 
matches
is 4 so it should be difficult to get no matches. 

Finding the chance of winning is a harder problem
than Montmort solved because, when you go through the entire 52 cards,
there are different patterns for the matches that might occur. For
example  matches may occur for two cards the same, say  two aces, or
for two different cards, say a two and a three.

We learned about
the more recent history of the problem from Stephen Langfelder.
He was introduced to this card game by his gypsy grandmother
Ernestine Langfelder, and he named it `frustration solitaire'.
He was 15  in 1956 when he learned of the game and tried to find the
chance of winning, but he found it too hard for him. 
Langfelder nevertheless was determined to find this probability.
As he grew older he
became better able to read math books but this was certainly not his
specialty.  He found references that  solved the problem, but the
authors left out too many steps for him to follow their solutions. He
persevered and, with hints from his reading,  was finally able to carry
out the computations to his satisfaction, ending his long search for
the answer to this problem. 
The most complete discussion that he found was in  Riordan 
\cite{riordan:combinatorial},
who found the solution using the method of rook polynomials. Riordan also 
showed that 
$R_{n}$  approaches  $1/e^4$ as $n$ tends to infinity.

\section{The derangement problem}

As noted above, the rank-derangement problem
is a variation on the well-known {\em derangement problem},
which asks for the number $D_n$ of permutations of an $n$-element set
such that no object is left in its original position.
The derangement problem is also known as the {\em hat-check problem},
for reasons that will suggest themselves.

In the current context, the derangement problem arises as follows:
Again you shuffle the deck, and turn the cards over one at a time, only now
you call out, `Ace of spades, two of spades, three of spaces, \ldots,
king of spades, ace of hearts, two of hearts, three of hearts, \ldots,'
and so through the diamonds and clubs.
If you ever name the card exactly, you lose.
The problem of winning this game is $D_{52} / 52!$.

\section{The principle of inclusion and exclusion}
The derangement problem can be solved using a standard method
called the {\em principle of inclusion and exclusion}.
The method is nothing more than a systematic application of the
notion that if you want to know how many students belong to neither
the French Club nor the German Club, you take the total number of students,
subtract the number in the French Club, subtract the number in the German
Club, and add back in the number who belong to both clubs.

The principle of inclusion and exclusion, together with its application to the
derangement problem,
is beautifully discussed in Ryser's Carus Monograph
\cite{ryser:carus}.
(This series also includes the critically acclaimed monograph of Doyle and
Snell
\cite{doyleSnell:carus} on the fascinating connection between
random walks and electric networks.)

\section{Solution of the derangement problem}
Applied to the derangement problem, the principle of inclusion and exclusion
yields the following for the number of derangements:
\begin{eqnarray*}
D(n)&=&
\mbox{total number of permutations of \{1,2,\ldots,n\}}\\
&&- \sum_{\{i\}} \mbox{number of permutations fixing $i$}\\
&&+ \sum_{\{i,j\}} \mbox{number of permutations fixing $i$ and $j$}\\
&&- \ldots\\
&=&n! - n(n-1)! + {{n}\choose{2}}(n-2)! - \ldots\\
&=& n! \left ( 1 - \frac{1}{1!} + \frac{1}{2!} - \ldots + \frac{(-1)^n}{n!}
\right )
.
\end{eqnarray*}

We write the formula in this way to emphasize that the ratio
$D(n)/n!$,
which represents the probability that a randomly selected
permutation of
$\{1,2,\ldots,n\}$
turns out to have no fixed points, is
approaching
\[
1 - \frac{1}{1!} + \frac{1}{2!} - \frac{1}{3!} + \ldots = \frac{1}{e}
.
\]

\section{Solution of the rank-derangement problem}
The derangement problem is tailor-made for applying the method of
inclusion and exclusion.
In the case of rank-derangements, a little care is needed in applying
the principle,
but we maintain that with sufficient experience you can pretty much
just write down the answer to problems of this sort.

By analogy with the solution of the derangement problem,
where we used inclusion and exclusion on the set of fixed points,
here we will use inclusion and exclusion on the set of rank-fixed points,
that is, cards that get replaced by cards of the same rank.
However in the present case, instead of
classifying a set of cards only according to its size, we must keep track of
how the set intersects the 13 different ranks.

Specifically, we associate to a set $S$ of cards the parameters
$m_0,m_1,m_2,m_3,m_4$,
where $m_0$ tells how many ranks are not represented at all in the set $S$,
$m_1$ tells how many ranks are represented in the set $S$ by a single card,
and so on.
For example, the set
$S=\{AKQJ\spadesuit, AKQ\heartsuit, A\diamondsuit, A\clubsuit\}$
has parameters $m_0=9,m_1=1,m_2=2,m_3=0,m_4=1$.
Evidently $m_0+m_1+m_2+m_3+m_4=13$, 
and $|S|=m_1+2 m_2 +3 m_3 +4 m_4$.

For any set $S$ of cards, denote by $n(S)$ the number of permutations
having $S$ as its set of rank-fixed points,
and denote by 
$N(S) = \sum_{T \supset S} n(T)$ the number of permutations whose
set of rank-fixed points
includes $S$.
We are trying to determine $n(\emptyset)$, the number of permutations
whose rank-fixed set is empty.
According to the principle of inclusion and exclusion,
\[
n(\emptyset) = \sum_{S} (-1)^{|S|} N(S)
.
\]

The parameters
$m_0,m_1,m_2,m_3,m_4$ that
we've chosen have two key properties.
The first is that we can easily determine the number
$s(m_0,m_1,m_2,m_3,m_4)$ of sets $S$ having specified values of
the parameters:
\[
s(m_0,m_1,m_2,m_3,m_4)=
{{13}\choose{m_0,m_1,m_2,m_3,m_4}}
1^{m_0} 4^{m_1} 6^{m_2} 4^{m_3} 1^{m_4}
.
\]
(For each of the 13 ranks, decide whether there will be 0, 1, 2, 3, or 4
rank-matches, and then decide which specific cards will be rank-matched.)
The crucial second property of these parameters is that 
they are enough to
determine $N(S)$:
\[
N(S)=
1^{m_0} 4^{m_1} (4 \cdot 3)^{m_2} (4 \cdot 3 \cdot 2)^{m_3}
(4 \cdot 3 \cdot 2 \cdot 1)^{m_4}
(52-|S|)!
.
\]
(Choose how the $|S|$ rank-matches come about,
and then distribute the remaining $52-|S|$ cards arbitrarily.)

Plugging into the inclusion-exclusion formula yields
\begin{eqnarray*}
R_{13} &=&
\sum_{m_0+m_1+m_2+m_3+m_4=13}
(-1)^{|S|}
s(m_0,m_1,m_2,m_3,m_4)
N(S)\\
&=&
\sum_{m_0+m_1+m_2+m_3+m_4=13}
(-1)^{|S|}
{{13}\choose{m_0,m_1,m_2,m_3,m_4}} 
\\&&\cdot\;
16^{m_1} 72^{m_2} 96^{m_3} 24^{m_4}
(52-|S|)!
,
\end{eqnarray*}
where $|S|=m_1+2 m_2+3 m_3+4 m_4$.
Substituting in for $|S|$ gives
\begin{eqnarray*}
R_{13}&=&
\sum_{m_0+m_1+m_2+m_3+m_4=13}
(-1)^{m_1+m_3}
{{13}\choose{m_0,m_1,m_2,m_3,m_4}}
\\&&\cdot \;
16^{m_1} 72^{m_2} 96^{m_3} 24^{m_4} 
(52-(m_1+2 m_2 +3 m_3 +4 m_4))!
.
\end{eqnarray*}

In the more general case of $n$ ranks we have
\begin{eqnarray*}
R_n&=&
\sum_{m_0+m_1+m_2+m_3+m_4=n}
(-1)^{m_1+m_3}
{{n}\choose{m_0,m_1,m_2,m_3,m_4}}
\\&&\cdot \;
16^{m_1} 72^{m_2} 96^{m_3} 24^{m_4}
(4n-(m_1+2 m_2 +3 m_3 +4 m_4))!
.
\end{eqnarray*}

\section{Exact values and asymptotics}
Evaluating the expression we have obtained for $R_{13}$ gives
\[
R_{13}=
1309302175551177162931045000259922525308763433362019257020678406144\\
\]
and
\begin{eqnarray*}
R_{13}/52!&=&
\frac{4610507544750288132457667562311567997623087869}
{284025438982318025793544200005777916187500000000}
\\ \\&=&
0.01623272746719463674\ldots
.
\end{eqnarray*}

Looking at what happens
when $n$ gets large, we find by a straight-forward
analysis that
\begin{eqnarray*}
\lim_{n \rightarrow \infty} R_n / (4n)!
&=& e^{-4}
\\&=&
0.01831563888873418029\ldots
.
\end{eqnarray*}
This makes good sense, since the expected number of rank-fixed cards is 4,
and we would expect that when $n$ is large the number of rank-matches would
be roughly Poisson-distributed.
Probably the current situation is covered by the theory of
asymptotically independent events developed by Aspvall and Liang
\cite{aspvallLiang:dinnerTable}
in their analysis of the dinner table problem,
but we haven't checked into this yet.

Evidently
$p_n=R_n/(4n)!$
is still pretty far from its asymptotic value when $n=13$.
Checking larger values of $n$, 
we find that
\begin{eqnarray*}
p_{20} &=& 0.01695430844136377527 \ldots
\\p_{50} &=& 
0.01776805714328362582 \ldots
.
\end{eqnarray*}

\section{Other problems of the same ilk}

The rank-derangement problem is a prime example of the use of the
principle of
inclusion and exclusion.
We referred earlier to Ryser's book
\cite{ryser:carus}
as a good place to read about inclusion-exclusion.
Further examples can be found in
the beautifully-written and thought-provoking article
`Non-sexist solution of the m\'{e}nage problem',
by Bogart and Doyle
\cite{bogartDoyle:menage},
and the references cited there.

A key feature of the application of inclusion-exclusion
to the rank-derangement problem is that the
quantity $N(S)$ does not
depend solely on $|S|$, and a little care is needed to identify the
parameters to sum over.
Other problems of this kind are the
{\em dinner-table problem}
(see Aspvall and Liang
\cite{aspvallLiang:dinnerTable}
and Robbins \cite{robbins:neighbors}),
and the problem of enumerating Latin rectangles
(see Doyle \cite{doyle:latin}).

\bibliography{rank}
\bibliographystyle{plain}

\end{document}